\documentclass{elsart}

\usepackage{amsmath}

\newcommand{\R}{\mathbf{R}}

\newcommand{\C}{\mathcal{C}}
\newcommand{\G}{\mathcal{G}}
\renewcommand{\L}{\mathcal{L}}
\DeclareMathOperator{\erf}{erf}

\renewcommand{\bar}{\tilde}

\begin{document}
\raggedbottom

\begin{frontmatter}
\title{Exact variations for stochastic heat equations driven by space--time 
white noise}
\author{Jan Pospisil}
\address{University of West Bohemia, New Technologies Research Centre, \\
306 14 Plzen, Czech Republic}
\ead{jan.pospisil@ntc.zcu.cz}
\author{Roger Tribe}
\address{University of Warwick, Mathematics Institute, \\
CV4 7AL Coventry, United Kingdom}
\ead{tribe@maths.warwick.ac.uk}
\thanks{This work was partially supported by the GACR Grant 201/04/0750 and by the 
MSMT research and development project LN00B084}
\begin{abstract}
This paper calculates the exact  
quadratic variation in space and quartic variation in time for the
solutions to a one dimensional stochastic heat equation
driven by a multiplicative space-time white noise. 
\end{abstract}
\end{frontmatter}

\section{Introduction}
We consider solutions $(X_t(x):t \geq 0, x \in \R)$ to 
the one-dimensional stochastic heat equation 
\begin{equation}\label{spde}
dX = \alpha \, \Delta X \, dt + \sigma(X) \, d \eta
\end{equation}
where $\alpha >0$, $\sigma:\R \to \R$ and $\eta(dx,dt)$ is a space-time white noise
on $\R \times [0,\infty)$. 
Such equations arise in several settings, for example generalised
Edwards--Wilkinson models for the roughening of surfaces \cite{Barabasi+Stanley}, 
continuum limits of particle processes (see for example \cite{Giacomin}, 
\cite{Konno+Shiga}, \cite{Tribe}), and 
continuous space parabolic Anderson models \cite{Andref}. 

Under reasonable conditions on $\sigma$, there exist continuous solutions 
which are known to be H\"{o}lder continuous of order $\frac12 - \epsilon$ 
in space and of order $\frac14 - \epsilon$ in time, for any $\epsilon>0$.  
The aim of this note is to calculate the exact variations.
We show that the paths $x \to X_t(x)$ have finite quadratic variation, for fixed $t>0$
and over an interval $[A_1,A_2]$, given by
\[
\lim_{n \to \infty} \sum_{j=1}^n \left( X_t(x_j) - X_t(x_{j-1}) \right)^2 = 
\frac{1}{2\alpha} \int^{A_2}_{A_1} \sigma^2(X_t(x)) dx,
\]
and the paths $t \to X_t(x)$ have quartic variation, over an interval $[T_1,T_2]$, given by 
\[
\lim_{n \to \infty} \sum_{j=1}^n \left( X_{t_j}(x) - X_{t_{j-1}}(x) \right)^4 = 
\frac{3}{\pi \alpha} \int^{T_2}_{T_1} \sigma^4(X_t(x)) dt.
\]
Here $\{x_j\}$ and $\{t_j\}$ are equally spaced partitions and the 
convergence is in probability. 
Precise statements are given later in Theorems \ref{T1} and \ref{T2}.
The results are local in character, being 
dependent only on the nature of the noise and the small time 
asymptotics of the Green's kernel, and we will remark later on similar results 
that hold for related stochastic PDEs.
In the case of constant $\sigma$ the variations can be checked by an $\L^2$ argument 
using the Gaussianity of the solution. For general $\sigma$ the proof involves 
showing that the increments are approximately conditionally Gaussian.

The result is interesting when considered as a method for parameter estimation since the value of the 
drift parameter $\alpha$ can theoretically be found from an arbitrarily small segment of a single path.
For stochastic ODEs the drift parameters are typically not reflected in the 
local structure of paths, since changes in value of the drift produce absolutely
continuous changes in the law of the process on finite time intervals. 
Here however, exploiting say the time variation, we may choose
\[
\hat{\alpha}_n = \frac{3(T_2-T_2)}{n\pi} \; \frac{\sum_{j=1}^n 
\sigma^4(X_{t_j}(x))}{\sum_{j=1}^n \left( X_{t_j}(x) - X_{t_{j-1}}(x) \right)^4}.
\]
as an estimator of $\alpha$, and an analogous estimator is possible 
that exploits the spatial variation. Consideration of the estimates involved in
our results suggests, at least when $\sigma$ is smooth, bounded and 
bounded away from zero, that 
$E[|\hat{\alpha}_n - \alpha| \wedge 1] = O(n^{-3/20})$.
Details on this and other parameter estimation problems will be found 
in the forthcoming thesis \cite{Pospisil}.  

\section{The linear case}
In this section we suppose that $\sigma(x)$ is a constant. By 
a linear change of variable we may suppose that $\sigma(x)=\alpha=1$.
The solution of (\ref{spde}) is then defined as
\begin{equation} \label{Green}
X_t(x) = \int_{\R} G_t(x-u) X_0(u) du \, + \int_0^t \int_{\R} G_{t-r}(x-u) \, \eta(du,dr) 
\end{equation}
where $ G_t(x) = \frac1{2 \sqrt{\pi t}} \exp (- x^2 / 4t)$
is the heat kernel. 
If we suppose some growth rate on the initial condition $X_0$,
then the first integral on the right hand side of (\ref{Green}) 
is a smooth function of $x \in \R$ and $t>0$ and 
will not affect the variations we consider. Without loss of generality we therefore assume 
for the rest of this section that $X_0(x)=0$.

The covariance structure of the solution (\ref{Green}), with $X_0=0$, 
can be calculated, using the isometry for space-time white noise integrals 
(see Walsh \cite{Walsh} chapter 2) and the 
semigroup property of the heat kernel, as, for $s \leq t$,
\begin{equation} \label{covariance}
E \left[ X_t(x)X_s(y) \right] = \int_0^s \frac{1}{2 \sqrt{\pi(s+t-2r)}} 
\exp \left(- \frac{(x-y)^2}{4(s+t-2r)} \right) dr.
\end{equation} 
\subsection{Spatial variation}
Throughout this subsection we fix $t>0$ and $A_1<A_2$.  
We consider the quadratic variation of $x \to X_t(x)$ 
over the interval $[A_1,A_2]$. From the covariance structure (\ref{covariance}) we 
may deduce, by integration by parts, 
\begin{multline} \label{quadraticcovariance}
E \left[X_t(y)X_t(x) \right] \\
= \sqrt{\frac{t}{2 \pi}} \exp \left( -\frac{(y-x)^2}{8t} \right) + 
\frac{(y-x)}{4} \erf \left( \frac{y-x}{2\sqrt{2t}} \right) - \frac{|y-x|}{4},
\end{multline}
where $\erf$ denotes the error function defined by 
$\erf(x) = (2/\sqrt{\pi}) \int_0^x e^{-t^2} dt$. 
In particular, taking $x=y$, we have $E[X^2_t(x)] = \sqrt{t/2\pi}$.
Define, for $\delta>0$, the jointly Gaussian increments
\[
\Delta(x,\delta) = X_t(x+ \delta)-X_t(x), \quad
\Delta(y,\delta) = X_t(y+\delta)-X_t(y).
\]
From the covariance structure (\ref{quadraticcovariance}) we have that
\begin{equation} \label{L11}
E \left[\Delta^2(x,\delta) \right]  = 
 \sqrt{\frac{2t}{\pi}} -
\sqrt{\frac{2t}{\pi}} \exp \left(-\frac{\delta^2}{8t} \right) - 
\frac{\delta}{2} \erf \left(\frac{\delta}{2\sqrt{2t}} \right) + \frac{\delta}{2}
= \frac{\delta}{2} + O(\delta^2).
\end{equation}
Here we have expanded in small $\delta$ and used the expansion
$ \erf(\delta) = \frac{2}{\sqrt{\pi}} \delta + O(\delta^3)$.
Furthermore, supposing that $x+\delta \leq y$, we have
\par

\begin{equation*}
\begin{split}
& E \left[ \Delta(x,\delta)\Delta(y,\delta) \right] \\
& = E \left[ X_t(x+\delta) X_t(y+\delta) - X_t(x+\delta)X_t(y) 
- X_t(x)X_t(y+\delta) + X_t(x)X_t(y) \right] \\
& = \sqrt{\frac{t}{2 \pi}} \left[ 2 \exp \left( -\frac{(y-x)^2}{8t} \right) 
  - \exp \left( -\frac{(y-x+\delta)^2}{8t} \right) - \right. \\
&\quad - \left. \exp\left( -\frac{(y-x-\delta)^2}{8t} \right) \right] 
  +  2 \, \frac{(y-x)}{4} \erf \left(\frac{y-x}{2\sqrt{2t}} \right) \\
&\quad - \frac{(y-x+\delta)}{4}\erf \left(\frac{y-x+\delta}{2\sqrt{2t}} \right) 
  - \frac{(y-x-\delta)}{4}\erf \left(\frac{y-x-\delta}{2\sqrt{2t}} \right).
\end{split}
\end{equation*}
Using the Taylor expansions
\begin{align*}
\exp \left( -\frac{(y-x \pm \delta)^2}{8t}\right) 
& = \exp \left(-\frac{(y-x)^2}{8t} \right) \left( 1 \mp 
   \frac{y-x}{4t} \, \delta  + O(\delta^2) \right) \\
\erf \left( \frac{y-x \pm \delta}{2\sqrt{2t}} \right) 
& = \erf \left(\frac{y-x}{2\sqrt{2t}} \right) \pm
   \frac{\delta}{\sqrt{2 \pi t}} \exp \left( - \frac{(y-x)^2}{8t} \right) 
   + O(\delta^2).
\end{align*}
we conclude that
\[
E \left[ \Delta(x,\delta) \Delta(y,\delta) \right] = O(\delta^{2}).
\]
Here we have expanded in small $\delta$ and used $t,t^{-1},x,y = O(1)$ where
appropriate. The expression $O(\delta^2)$ denotes a quantity bounded by
$C \delta^2$, for small $\delta$, with the same constant $C$ for all 
$x,y \in [A_1,A_2]$ satisfying $x+\delta \leq y$. 
We will need two higher Gaussian moments which can be read off from the lower moments
(for example using $ E[X^2 Y^2] = E[X^2] \, E[Y^2]  + 2(E[XY])^2$). When
$x+\delta \leq y$ we have
\begin{align} 
E \left[\Delta^4(x,\delta) \right] & = O(\delta^2), \label{L12} \\ 
E \left[ \Delta^2(x,\delta) \Delta^2(y,\delta) \right] 
& = \frac{\delta^2}{4} + O(\delta^{3}). \label{L13}
\end{align}
\begin{prop} \label{linearquadratic}
Define, for $j=0,1,\ldots,n$, a space grid by 
$ x_j = A_1+ j \delta$, where $\delta= \frac{1}{n} (A_2-A_1)$. 
Then the following limit holds in mean square:
\[
\lim_{n \to \infty} \; \sum_{j=1}^{n} \left(X_t(x_j)-X_t(x_{j-1}) \right)^2 = 
\frac12 (A_2-A_1).
\]
\end{prop}

\begin{pf}
The desired $\L^2$ convergence can be rewritten as
\[
\lim_{n \to \infty} \sum_{j=1}^{n} \sum_{k=1}^{n} 
E \left[  \left( \Delta^2(x_j,\delta)-\frac{\delta}{2} \right) 
\left(\Delta^2(x_k,\delta)-\frac{\delta}{2}\right) \right] =0.
\]
Consider first the diagonal terms when $j=k$. These give a contribution bounded by
\[
\sum_{j=1}^{n} E\left[ \Delta^4(x_j,\delta) + \delta \Delta^2(x_j,\delta) 
+ \frac{\delta^2}{4} \right] 
\leq n O(\delta^2) \to 0 \quad \mbox{as $n \to \infty$}
\]
using the estimates from (\ref{L11},\ref{L12}) and the fact that $\delta = O(1/n)$.
The off-diagonal terms give
\begin{multline*}
2 \sum_{j=1}^{n}\sum_{k=j+1}^{n} 
E \left[ \Delta^2(x_j,\delta) \Delta^2(x_k,\delta) - 
\frac{\delta}{2}(\Delta^2(x_j,\delta) + \Delta^2(x_k,\delta)) + \frac{\delta^2}{4} \right] \\
= 2 \sum_{j=1}^{n}\sum_{k=j+1}^{n} O(\delta^{3}).
\end{multline*}
The cancellation in the last equality uses the expansions in (\ref{L11},\ref{L13}).
The final sum converges to zero as $n \to \infty$ completing the proof.\qed
\end{pf}
\subsection{Temporal variation}
Throughout this subsection 
we fix $x \in \R$ and times $0<T_1<T_2$. We consider the quartic variation of $t \to X_t(x)$ 
over the interval $[T_1,T_2]$. From the covariance structure (\ref{covariance}) we deduce,
for $s<t$, 
\[
E \left[X_t(x)X_s(x) \right]  =  \frac1{2\sqrt{\pi}} \left(\sqrt{t+s}-\sqrt{t-s}\right).
\]
Suppose that $s,t \in [T_1,T_2]$ and $0 < \delta \leq t-s$.
Define the jointly Gaussian increments
\[
\Delta(s,\delta) = X_{s+\delta}(x)-X_s(x), \quad
\Delta(t,\delta) = X_{t+\delta}(x)-X_t(x).
\]
Using the covariance structure above we find 
\[
E \left[\Delta^2(s,\delta) \right]  = 
\frac1{2 \sqrt{\pi}} \left( \sqrt{2(s+\delta)}+\sqrt{2s}-2\sqrt{2s+\delta}+2\sqrt{\delta} \right) 
= \frac{\delta^{1/2}}{\sqrt{\pi}} + O(\delta^2).
\]
Furthermore
\begin{equation*}
\begin{split}
E \left[ \Delta(s,\delta)\Delta(t,\delta) \right] & = 
  -\frac1{2\sqrt{\pi}} \left[ 2 \sqrt{t+s+\delta} - \sqrt{t+s+2\delta} - \sqrt{t+s} \right. \\
&\quad + \left. 2\sqrt{t-s} - \sqrt{t-s+\delta} - \sqrt{t-s-\delta} \right] \\
& = O( \delta^{2}) + O((t-s)^{-3/2} \delta^2 ).
\end{split}
\end{equation*}
Here we have expanded as a series in $\delta$ and $(t-s)^{-1}$ using that 
$s,t,s+t = O(1)$ where appropriate. We deduce from this covariance structure
the following higher asymptotics, which hold uniformly over
$s,t \in [T_1,T_2]$ satisfying $\delta \leq t-s$,
\begin{align}
E \left[ \Delta^4(s,\delta) \right] & =  
 \frac{3 \delta}{\pi} + O(\delta^{5/2}), \label{L21} \\
E \left[ \Delta^8(s,\delta) \right] & =  O(\delta^{2}), \label{L22} \\ 
E \left[ \Delta^4(s,\delta) \Delta^4 (t,\delta) \right] & =  
\frac{9 \delta^2}{\pi^2} + O(\delta^{7/2}) + O((t-s)^{-3} \delta^5). \label{L23}
\end{align}
We have used the expansion of Gaussian higher moments in terms of the
covariance function, for example
$ E[X^4 Y^4] = 9( E[X^2] E[Y^2])^2 + 24(E[XY])^4 + 72 E[X^2] E[Y^2](E[XY])^2 $,
keeping only the leading terms of the asymptotics in $\delta$ and $(t-s)^{-1}$ that
we will need.
\begin{prop} \label{linearquartic}
Define, for $j=0,1,\ldots,n$, a time grid by 
$ t_j = T_1+j \delta$, where $\delta= \frac{1}{n} (T_2-T_1)$. 
Then the following limit holds in mean square:
\[
\lim_{n \to\infty} \; \sum_{j=1}^{n} \left(X_{t_j}(x)-X_{t_{j-1}}(x) \right)^4 = 
\frac{3}{\pi} (T_2-T_1).
\]
\end{prop}
\begin{pf} 
The desired $\L^2$ convergence can be rewritten as
\[
\lim_{n \to \infty} \sum_{j=1}^{n} \sum_{k=1}^{n} 
E \left[  \left( \Delta^4(t_j,\delta)-\frac{3 \delta}{\pi} \right) 
\left(\Delta^4(t_k,\delta)-\frac{3 \delta}{\pi}\right) \right] =0.
\]
Consider first the diagonal terms when $j=k$. These give a contribution bounded by
\[
\sum_{j=1}^{n} E\left[ \Delta^8(t_j,\delta) + \frac{6 \delta}{\pi} \Delta^4(t_j,\delta) 
+ \frac{9 \delta^2}{\pi^2} \right] 
\leq n O(\delta^2) \to 0 \quad \mbox{as $n \to \infty$}
\]
using the estimates from (\ref{L21},\ref{L22}) and the fact that $\delta = O(1/n)$.
The off-diagonal terms give
\begin{multline*}
2 \sum_{j=1}^{n}\sum_{k=j+1}^{n} 
E\left[ \Delta^4(t_j,\delta) \Delta^4(t_k,\delta) 
- \frac{3 \delta}{\pi}(\Delta^4(t_j,\delta) + \Delta^4(t_k,\delta)) 
+ \frac{9 \delta^2}{\pi^2} \right] \\
= 2 \sum_{j=1}^{n}\sum_{k=j+1}^{n} \left( 
 (k-j)^{-3} O(\delta^2) + O(\delta^{7/2}) \right). \\
\end{multline*}
The cancellation of the moments in the last equality uses the expansions from 
(\ref{L21},\ref{L23}).
The final sum converges to zero as $n \to \infty$ completing the proof.\qed
\end{pf} 

\section{The non-linear case} 
We consider the equation (\ref{spde}) with $\sigma$ globally Lipschitz and satisfying a 
linear growth condition. By linear time scaling we may assume again that $\alpha=1$.
We suppose that $X_0$ is in the following 
space $\C_{tem}$ of functions with slower than exponential growth
\[
\C_{tem} = \{f:\R \to \R: \|f\|_{\lambda} := \sup_x |f(x)| e^{-\lambda x} < \infty, \;
\mbox{for all $\lambda >0$} \}.
\] 
Give $\C_{tem}$ the topology induced by the seminorms $\|f\|_{\lambda}$ for $\lambda >0$. 
The choice of $\C_{tem}$ is not crucial but allows us to quote the following 
results.
Let $\eta(dx,dt)$ be a space-time white noise on $\R \times [0,\infty)$,
adapted to a filtered probability space $(\Omega, \mathcal{F}, 
(\mathcal{F})_t,P)$. 
In Shiga \cite{Shiga} Theorem 2.1 it is shown there exists
a pathwise unique adapted, continuous $\C_{tem}$ valued solution, which moreover satsifies the 
moment bounds, for any $p, \lambda, T >0$,
\begin{equation} \label{momentbounds}
E \left[ |X_t(x)|^p \right] \leq C(p, \lambda, T) e^{\lambda |x|} \quad 
\mbox{for all $x\in \R$ and $t \in [0,T]$,}
\end{equation}
and the increment estimates
\begin{equation} \label{incrementbounds}
E \left[ |X_t(x)-X_s(y)|^p \right] \leq C(p,\lambda,T)
\left( |x-y|^{p/2} + |t-s|^{p/4} \right) e^{\lambda (|x| + |y|)} 
\end{equation}
for all $x,y \in \R$ and $s,t \in [0,T]$.
In fact these, more or less standard, bounds are not explicitly stated in \cite{Shiga} 
but are implicit in the proofs found in the appendix. 

The solution can be written as $X_t(x) = \int G_t(x-u) X_0(u) du + \bar{X}_t(x)$ where
\begin{equation} \label{barXdefn}
\bar{X}_t(x) = \int_0^t \int_{\R} G_{t-r}(x-u) \sigma (X_r(u)) \, \eta(du, dr).
\end{equation}
Since $(t,x) \to \int G_t(x-u) X_0(u) du$ is smooth on $t>0$, 
the quadratic and quartic 
variations of $X_t(x)$ will exactly coincide with those of $\tilde{X}_t(x)$.
We therefore now restrict our consideration to the process $\bar{X}_t(x)$. 
We will calculate the spatial and temporal variation in subsections
\ref{s3.1} and \ref{s3.2} respectively. In \ref{s3.3} we describe
some methods to transfer these results to related equations. 
\subsection{Spatial variation} \label{s3.1}
We fix $t>0$ and consider the spatial increments defined by
\[
\Delta(X,x,\delta) = \bar{X}_t(x+\delta) - \bar{X}_t(x).
\]
The key idea is that the main contribution to $\Delta(X,x,\delta)$
comes from the noise near time $t$. We will approximate the increment 
$\Delta(X,x,\delta)$ by the term $\sigma(X_{t(\delta)}(x)) \bar{\Delta}(x,\delta)$ where 
$t(\delta) = t - \delta^{4/3}$ and
\begin{multline*}
\bar{\Delta}(x,\delta) = \int^t_{t(\delta)} \int 
\left( G_{t-r}(x+\delta-u) - G_{t-r}(x-u) \right) \eta(du,dr) \\
+ \int^{t(\delta)}_0 \int 
\left( G_{t-r}(x+\delta-u) - G_{t-r}(x-u) \right) \bar{\eta}(du,dr)
\end{multline*}
for $\bar{\eta}(dx,dt)$ an independent space-time white noise.
The integral against $\bar{\eta}(dx,dt)$ is included so that
$\bar{\Delta}(x,\delta)$ has the same Gaussian distribution 
as $\Delta(x,\delta)$ from the linear theory. Moreover 
$\bar{\Delta}(x,\delta)$ is independent of $\sigma\{X_s:s \leq t(\delta)\}$.
The value of $t(\delta)$ is chosen to optimise the estimate in the following lemma, which 
shows the above approximation is valid.
\begin{lem} \label{keyspacelemma}
For any $\lambda,T>0$ there exists $C(\lambda, T) < \infty$ so that
\[
E \left[ \left| \Delta(X,x,\delta) - 
\sigma(X_{t(\delta)}(x)) \bar{\Delta}(x,\delta) \right|^2 \right] 
\leq C(\lambda,T) \, \delta^{4/3} \, e^{\lambda |x|}
\]
for all $x \in \R$, $0 \leq \delta^{4/3} \leq t \wedge 1$ and $t \leq T$. 
\end{lem}
\begin{pf}
We split the difference into three parts: 
\begin{align*}
& \Delta(X,x,\delta) - \sigma(X_{t(\delta)}(x)) \bar{\Delta}(x,\delta) \\
& = \int\limits^t_{t(\delta)} \int \left( G_{t-r}(x+\delta-u) - G_{t-r}(x-u) \right) 
\left( \sigma(X_r(u)) - \sigma(X_{t(\delta)}(x)) \right) \eta(du,dr) \\
& \quad\quad+ \int\limits^{t(\delta)}_0 \int \left( G_{t-r}(x+\delta-u) - G_{t-r}(x-u) \right) 
\sigma(X_{r}(u)) \, \eta(du,dr) \\
& \quad\quad\quad\quad - \int\limits^{t(\delta)}_0 \int \left( G_{t-r}(x+\delta-u) - G_{t-r}(x-u) \right) 
\sigma(X_{t(\delta)}(x)) \, \bar{\eta}(du,dr) \\
& = I + II - III.
\end{align*}
We estimate each of these terms separately.
We need the deterministic estimates, for $\lambda>0$,
\begin{equation}
\int^t_s \int G^2_{t-r}(x-u) e^{\lambda |u|} du \, dr  \leq   
C(\lambda,T) (t-s)^{1/2} e^{\lambda |x|} \label{estimate1} \\
\end{equation}
and
\begin{multline}
\int^s_0 \int \left( G_{t-r}(x+\delta-u) - G_{t-r}(x-u) \right)^2 e^{\lambda |u|} du \, dr \\
 \leq  C(\lambda,T) \delta^2 (t-s)^{-1/2} e^{\lambda |x|}, \label{estimate2}
\end{multline}
valid whenever $\delta \in [0,1]$, $ 0 \leq s \leq t \leq T$ and  $ x \in \R$.
These (and the similar estimate (\ref{estimate3}) below) follow from straightforward 
calculations and we omit the proofs. Using the moments (\ref{momentbounds}) 
and estimate (\ref{estimate2}) we have
\begin{multline*}
E \left[ |II|^2 \right] =  
\int^{t(\delta)}_0 \int \left( G_{t-r}(x+\delta-u) - G_{t-r}(x-u) \right)^2 
E \left[ \sigma^2(X_r(u)) \right] du \, dr \\
\leq C(\lambda,T) e^{\lambda |x|} \delta^2 \left(t-t(\delta) \right)^{-1/2} = 
C(\lambda,T) e^{\lambda |x|} \delta^{4/3}. 
\end{multline*}
The estimate on $E[|III|^2]$ is entirely similar. Using the increment bounds
(\ref{incrementbounds}) we bound $E[|I|^2]$ by 
\begin{equation*}
\begin{split}
& C(\lambda,T) e^{\lambda |x|} \int\limits^t_{t(\delta)} \int 
\left( G_{t-r}(x+\delta-u) - G_{t-r}(x-u) \right)^2 
\left( |x-u| + \delta^{2/3} \right) e^{\lambda |u|} du \, dr \\
& \leq C(\lambda,T) e^{\lambda |x|} \int^t_{t(\delta)} \int 
 G^2_{t-r}(x-u) \left( |x-u| + \delta + \delta^{2/3} \right) e^{\lambda |u|} du \, dr\\
& \leq C(\lambda,T) e^{2 \lambda |x|} \delta^{4/3} + C(\lambda,T) e^{\lambda |x|} 
\int^t_{t(\delta)} \int G^2_{t-r}(x-u) |x-u|e^{\lambda |u|} du \, dr\\
& \leq C(\lambda,T) e^{2 \lambda |x|} \delta^{4/3} + C(\lambda,T) e^{\lambda |x|} 
\int^t_{t(\delta)} \int G_{t-r}(x-u) e^{\lambda |u|} du \, dr\\
& \leq C(\lambda,T) e^{2 \lambda |x|} \delta^{4/3}.
\end{split}
\end{equation*}
The second inequality here uses (\ref{estimate1}) and in the third we used
the bound $\sup_{t,z} |zG_t(z)| < \infty$.
Noting $\lambda>0$ is arbitrary, we see that 
combining the bounds on the terms
$I$, $II$ and $III$ completes the proof.\qed
\end{pf}

Using this approximation we can establish the spatial variation in the 
non-linear case along the same lines as for the linear case.
\begin{thm} \label{T1}
Fix $A_1 < A_2$ and $t>0$. Define, for $j=0,1,\ldots,n$, a space grid by 
$ x_j = A_1+ j \delta$, where $\delta= \frac{1}{n} (A_2-A_1)$. 
Then the following limit holds in probability:
\[
\lim_{n \to \infty} \; \sum_{j=1}^{n} \left(\bar{X}_t(x_j) - 
\bar{X}_t(x_{j-1}) \right)^2 = \frac12 \int^{A_2}_{A_1} \sigma^2(X_t(u)) du.
\]
\end{thm}
\begin{pf}
We break the required convergence into three parts as follows.
\begin{align}
& \sum_{j=1}^n \Delta^2(X,x_j,\delta) 
- \frac12 \int^{A_2}_{A_1} \sigma^2(X_t(u)) du \nonumber \\
& \quad = \sum_{j=1}^n \left( \Delta^2(X,x_j,\delta) 
- \sigma^2(X_{t(\delta)}(x_j)) \bar{\Delta}^2(x_j,\delta) \right) \label{term1s} \\
& \quad\quad\quad + \sum_{j=1}^n \sigma^2(X_{t(\delta)}(x_j)) \left( \bar{\Delta}^2(x_j,\delta) 
- \frac{\delta}{2} \right) \label{term2s} \\
& \quad\quad\quad\quad\quad + \frac{\delta}{2} \sum_{j=1}^n  \sigma^2(X_{t(\delta)}(x_j)) 
- \frac12 \int^{A_2}_{A_1} \sigma^2(X_t(u)) du. \label{term3s} 
\end{align}
The third term (\ref{term3s}) converges almost surely to zero as $n \to \infty$
using the uniform continuity of $(s,x) \to X_s(x)$ and
a Riemann sum approximation to the integral.
The first term (\ref{term1s}) converges to zero in $\L^1(P)$ since
\begin{multline*}
\sum_{j=1}^n E \left[ \left| \Delta^2(X,x_j,\delta) 
- \sigma^2(X_{t(\delta)}(x_j)) \bar{\Delta}^2(x_j,\delta) \right| \right] \\
\begin{aligned}
& \leq \sum_{j=1}^n \left( E \left[ | \Delta(X,x_j,\delta) 
- \sigma(X_{t(\delta)}(x_j)) \bar{\Delta}(x_j,\delta)|^2 \right] \right)^{1/2} \\
&\quad\quad\quad\cdot\left( E \left[ | \Delta(X,x_j,\delta) 
+ \sigma(X_{t(\delta)}(x_j)) \bar{\Delta}(x_j,\delta) |^2 \right] \right)^{1/2} \\
& \leq C(t,A_1,A_2) \delta^{2/3} \sum_{j=1}^n 
\left( E \left[  \Delta^2(X,x_j,\delta) 
+ \sigma^2(X_{t(\delta)}(x_j)) \bar{\Delta}^2(x_j,\delta) \right] \right)^{1/2} \\
& \leq C(t,A_1,A_2) \delta^{2/3}  n O (\delta^{1/2}) \to 0.
\end{aligned}
\end{multline*}
The first inequality is Cauchy-Schwartz, 
the second inequality uses the bound from Lemma \ref{keyspacelemma}
and in the final inequality we have applied the moments 
from (\ref{momentbounds}) and (\ref{incrementbounds}) and the independence between 
$\bar{\Delta}(x_j,\delta)$ and $\sigma(X_{t(\delta)}(x_j))$.

We will now show the second term (\ref{term2s}) converges to zero in
$\L^2(P)$ by mimicking the linear case. Indeed, using the independence between 
$\bar{\Delta}(x_j,\delta)$ and $\sigma(X_{t(\delta)}(x_j))$ and the 
moments (\ref{L11},\ref{L12}), we have
\[
E \left[\sum_{j=1}^n  \sigma^4(X_{t(\delta)}(x_j)) \left( \bar{\Delta}^2(x_j,\delta) 
- \frac{\delta}{2} \right)^2 \right]
 \leq  E \left[\sum_{j=1}^n  \sigma^4(X_{t(\delta)}(x_j)) O(\delta^2) \right] 
\]
which converges to zero using the moments bounds (\ref{momentbounds}). Similarly, using
(\ref{L13}),
\begin{multline*}
E \left[\sum_{j,k=1, k \neq j}^n  
\sigma^2(X_{t(\delta)}(x_j)) \sigma^2(X_{t(\delta)}(x_k))  
\left( \bar{\Delta}^2(x_j,\delta) - \frac{\delta}{2} \right)
\left( \bar{\Delta}^2(x_k,\delta) - \frac{\delta}{2} \right)
 \right] \\ 
= E \left[\sum_{j,k=1, k \neq j}^n  
\sigma^2(X_{t(\delta)}(x_j)) \sigma^2(X_{t(\delta)}(x_k))  O(\delta^{3})  \right]
\end{multline*}
which also converges to zero. This shows the second term (\ref{term2s})
converges to zero in $\L^2(P)$ and completes the proof.\qed
\end{pf}
\subsection{Temporal variation} \label{s3.2}
Fix $x$ and define temporal increments via
\[
\Delta(X,t,\delta) = \bar{X}_{t+\delta}(x) - \bar{X}_t(x).
\]
The key idea is that the main contribution to $\Delta(X,t,\delta)$
comes from the noise near time $t$. We will approximate the increment 
by $\sigma(X_{t(\delta)}(x)) \bar{\Delta}(t,\delta)$ where 
$t(\delta) = t - \delta^{4/5}$ and 
\begin{multline} \label{timeapprox}
\bar{\Delta}(t,\delta) = \int^{t+\delta}_{t(\delta)} \int 
\left( G_{t+\delta-r}(x-u) - G_{t-r}(x-u) \right) \eta(du,dr) \\
+ \int^{t(\delta)}_0 \int 
\left( G_{t+\delta-r}(x-u) - G_{t-r}(x-u) \right) \bar{\eta}(du,dr). 
\end{multline}
where $\bar{\eta}(dx,dt)$ is an independent space-time white noise.
We use here, and below, the convention that $G_t(x)=0$ for $t<0$.
Note again that $\bar{\Delta}(t,\delta)$ has the same Gaussian distribution 
as $\Delta(t,\delta)$ from the linear theory and 
is independent of $\sigma\{X_s:s \leq t(\delta)\}$.
The value of $t(\delta)$ is chosen to optimise the estimate in the following lemma, which 
shows the above approximation is valid.
\begin{lem} \label{keytimelemma}
For any $\lambda,T>0$ there exists $C(\lambda, T) < \infty$ so that
\[
E \left[ \left| \Delta(X,t,\delta) - 
\sigma(X_{t(\delta)}(x)) \bar{\Delta}(t,\delta) \right|^2 \right] 
\leq C(\lambda,T) \delta^{4/5} e^{\lambda |x|}
\]
for all $x \in \R$, $0 \leq \delta^{4/5} \leq t \wedge 1$ and $t \leq T$. 
\end{lem}
\begin{pf}
We split the difference into three parts:  
\begin{multline*}
\Delta(X,t,\delta) - \sigma(X_{t(\delta)}(x)) \bar{\Delta}(t,\delta) \\
\begin{aligned}
& = \int\limits^{t+\delta}_{t(\delta)} \int \left( G_{t+\delta-r}(x-u) - G_{t-r}(x-u) \right) 
\left( \sigma(X_r(u)) - \sigma(X_{t(\delta)}(x)) \right) \eta(du,dr) \\
& \quad\quad + \int\limits^{t(\delta)}_0 \int \left( G_{t+\delta-r}(x-u) - G_{t-r}(x-u) \right) 
\sigma(X_{r}(u)) \, \eta(du,dr) \\
& \quad\quad\quad\quad - \int\limits^{t(\delta)}_0 \int \left( G_{t+\delta-r}(x-u) - G_{t-r}(x-u) \right) 
\sigma(X_{t(\delta)}(x)) \, \bar{\eta}(du,dr)\\
&= I + II - III.
\end{aligned}
\end{multline*}
We estimate each of these terms separately. 
We need the deterministic estimate, for $\lambda>0$,
\begin{multline}
\int^s_0 \int \left( G_{t+\delta-r}(x-u) - G_{t-r}(x-u) \right)^2 e^{\lambda |u|} du \, dr \\
\leq C(\lambda,T) \delta^2 (t-s)^{-3/2} e^{\lambda |x|}, \label{estimate3}
\end{multline}
valid whenever $\delta \in [0,1]$, $ 0 \leq s \leq t \leq T$ and $x \in \R$.
Using the moments (\ref{momentbounds}) and estimate (\ref{estimate3}) we have
\begin{multline*}
E \left[ |II|^2 \right] = 
\int^{t(\delta)}_0 \int \left( G_{t+\delta-r}(x-u) - G_{t-r}(x-u) \right)^2 
E \left[ \sigma^2(X_r(u)) \right] du \, dr \\
\leq C(\lambda,T) e^{\lambda |x|} \delta^2 \left(t-t(\delta) \right)^{-3/2} =
C(\lambda,T) e^{\lambda |x|} \delta^{4/5}.
\end{multline*}
The estimate on $E[|III|^2]$ is entirely similar. Using the increment bounds
(\ref{incrementbounds}) we bound $E \left[ |I|^2 \right]$ by 
\[ 
C(\lambda,T) e^{\lambda |x|} \int^{t+\delta}_{t(\delta)} \int 
\left( G_{t+\delta-r}(x-u) - G_{t-r}(x-u) \right)^2 
\left( |x-u| + \delta^{2/5} \right) e^{\lambda |u|} du \, dr 
\]
which, as arguing as in Lemma \ref{keyspacelemma}, is bounded by $C(\lambda,T) e^{2 \lambda |x|}
\delta^{4/5}$. Combining the bounds on the terms
$I$, $II$ and $III$ completes the proof.\qed
\end{pf}

Using this approximation we can establish the spatial variation in the 
non-linear case along the same lines as for the linear case.
\begin{thm}  \label{T2}
Fix $0<T_1 < T_2$ and $x \in \R$. Define, for $j=0,1,\ldots,n$, a time grid by 
$ t_j = T_1+ j \delta$, where $\delta= \frac{1}{n} (T_2-T_1)$. 
Then the following limit holds in probability:
\[
\lim_{n \to \infty} \; \sum_{j=1}^{n} \left(\bar{X}_{t_j}(x) - 
\bar{X}_{t_{j-1}}(x) \right)^4 = \frac{3}{\pi} \int^{T_2}_{T_1} \sigma^4(X_t(x)) dt.
\]
\end{thm}

\begin{pf}
We use one more trick in the approximation for $\Delta(X,t_j,\delta)$. 
We approximate this increment by $\sigma(X_{t_j(\delta)}(x)) \bar{\Delta}_j(t_j,x)$
where $t_j(\delta) = t_j - \delta^{4/5}$ and where $\bar{\Delta}_j(t_j,x)$ are 
defined as in (\ref{timeapprox}) using an I.I.D. 
sequence $(\bar{\eta}_j(dx,dt):j=0,1,\ldots)$ of independent space-time white noises. 

We break the required convergence into three parts as follows.
\begin{align}
& \sum_{j=1}^n \Delta^4(X,t_j,\delta) 
- \frac{3}{\pi} \int^{T_2}_{T_1} \sigma^4(X_t(x)) dt \nonumber \\
& \quad = \sum_{j=1}^n \left( \Delta^4(X,t_j,\delta) 
- \sigma^4(X_{t_j(\delta)}(x)) \bar{\Delta}_j^4(t_j,\delta) \right) \label{term1} \\
& \quad\quad\quad + \sum_{j=1}^n \sigma^4(X_{t_j(\delta)}(x)) \left( \bar{\Delta}_j^4(t_j,\delta) 
- \frac{3\delta}{\pi} \right) \label{term2} \\
& \quad\quad\quad\quad\quad + \frac{3\delta}{\pi} \sum_{j=1}^n  \sigma^4(X_{t_j(\delta)}(x)) 
- \frac{3}{\pi} \int^{T_2}_{T_1} \sigma^4(X_t(x)) dt. \label{term3} 
\end{align}
The third term (\ref{term3}) converges almost surely to zero as $n \to \infty$
using the uniform continuity of $(s,x) \to X_s(x)$ and
a Riemann sum approximation to the integral.
The first term (\ref{term1}) converges to zero in $\L(P)^1$ since
\begin{multline*}
\sum_{j=1}^n E \left[ \left| \Delta^4(X,t_j,\delta) 
- \sigma^4(X_{t_j(\delta)}(x)) \bar{\Delta}_j^4(t_j,\delta) \right| \right] \\
\begin{aligned}
& \leq C \sum_{j=1}^n \left( E \left[ | \Delta(X,t_j,\delta) 
- \sigma(X_{t_j(\delta)}(x)) \bar{\Delta}_j(t_j,\delta)|^2 \right] \right)^{1/2} \\
& \quad\quad\quad\cdot\left( E \left[ \Delta^6(X,t_j,\delta) 
+ \sigma^6(X_{t(\delta)}(t_j)) \bar{\Delta}_j^6(t_j,\delta) \right] \right)^{1/2} \\
& \leq C(|x|,T_2) \delta^{2/5} n \, O (\delta^{3/4}) \to 0. 
\end{aligned}
\end{multline*}
The final inequality uses the bound from Lemma \ref{keyspacelemma}, 
the increment moments  from (\ref{momentbounds}) and (\ref{incrementbounds}) and 
the independence between 
$\bar{\Delta}(t_j,\delta)$ and $\sigma(X_{t_j(\delta)}(x))$. 

We will now show the second term (\ref{term2}) converges to zero in
$\L^2(P)$ by mimicking the linear case. We consider the double sum
\begin{equation} \label{L2}
\sum_{j,k=1}^n E \left[ 
\sigma^4(X_{t_j(\delta)}(x)) \sigma^4(X_{t_k(\delta)}(x))
 \left( \bar{\Delta}_j^4(t_j,\delta) - \frac{3\delta}{\pi} \right)
 \left( \bar{\Delta}_k^4(t_k,\delta) - \frac{3\delta}{\pi} \right)
\right]. 
\end{equation}
For a term corresponding to $j,k$ satisfying  $|j-k| \leq n^{1/2} $ we 
apply H\"{o}lder and use the independence between 
$\bar{\Delta}(t_j,\delta)$ and $\sigma(X_{t_j(\delta)}(x))$
and the asymptotics (\ref{L21},\ref{L22}) for moments of $\Delta(t,\delta)$ from 
the linear case to find
\begin{multline*}
E \left[ 
\sigma^4(X_{t_j(\delta)}(x)) \sigma^4(X_{t_k(\delta)}(x))
 \left( \bar{\Delta}_j^4(t_j,\delta) - \frac{3\delta}{\pi} \right)
 \left( \bar{\Delta}_k^4(t_k,\delta) - \frac{3\delta}{\pi} \right)
\right] \\ 
\begin{aligned}
& \leq \left( E \left[ \sigma^8(X_{t_j(\delta)}(x)) 
\left( \bar{\Delta}_j^4(t_j,\delta) - \frac{3\delta}{\pi} \right)^2
\right] \right)^{1/2} \\
&\quad\cdot \left( E \left[ \sigma^8(X_{t_k(\delta)}(x))
 \left( \bar{\Delta}_k^4(t_k,\delta) - \frac{3\delta}{\pi} \right)^2
\right] \right)^{1/2} \\ 
& = O(\delta^2).
\end{aligned}
\end{multline*}
This implies that the sum over $j,k$ satisfying $|j-k| \leq n^{1/2} $ 
in (\ref{L2}) converges to zero as $n \to \infty$. 
Let $\G(k,j) = \sigma\{ X_s(x): s \leq t_k(\delta), x \in \R\} \vee
\sigma\{\bar{\eta}_j(dx,\,dt): t \geq 0, x \in \R\}$. For
$k \geq j + n^{1/2}$, the variable $\bar{\Delta}_k(t_k,\delta)$ is independent
of $\G(k,j)$ whereas $\bar{\Delta}_j(t_j,\delta)$ is $\G(k,j)$ 
measurable. For such $j,k$ we therefore have
\begin{multline*}
E \left[ 
\sigma^4(X_{t_j(\delta)}(x)) \sigma^4(X_{t_k(\delta)}(x))
 \left( \bar{\Delta}_j^4(t_j,\delta) - \frac{3\delta}{\pi} \right)
 \left( \bar{\Delta}_k^4(t_k,\delta) - \frac{3\delta}{\pi} \right)
\right] \\ 
\begin{aligned}
& \leq E \left[ 
\sigma^4(X_{t_j(\delta)}(x)) \sigma^4(X_{t_k(\delta)}(x))
 \left( \bar{\Delta}^4(t_j,\delta) - \frac{3\delta}{\pi} \right)
O(\delta^{5/2}) \right] \\ 
& \leq O(\delta^{5/2}) \left( E \left[ 
\sigma^8(X_{t_j(\delta)}(x)) \sigma^8(X_{t_k(\delta)}(x)) \right] \right)^{1/2}
\left( E \left[ \left( \bar{\Delta}_j^4(t_j,\delta) - \frac{3\delta}{\pi} \right)^2
\right] \right)^{1/2} \\
& = O(\delta^{9/2}). 
\end{aligned}
\end{multline*}
For the first inequality here we conditioned
on $\G(k,j)$ and used the asymptotics from the linear case to calculate
$E[ \bar{\Delta}_k^4(t_k,x) - (3 \delta/\pi) ]$.
This implies that the sum over $j,k$ satisfying $k \geq j + n^{1/2}$ 
in (\ref{L2}) converges to zero as $n \to \infty$. 
A similar estimate holds for the terms $ j  \geq k + n^{1/2} $
and this completes the proof.\qed
\end{pf}
\subsection{Related results} \label{s3.3}
The results above rely only on the correlation of the space-time white noise
and the small $t$ asymptotics of the heat kernel $G_t(x)$, and so hold for
many variations of the stochastic PDE considered above. We comment here on some 
possible changes. 

The addition of a drift term such as 
\begin{equation} \label{withdrift}
dX = \alpha \, \Delta X \, dt + b(X)\,dt + \sigma(X) \, d \eta
\end{equation}
will not affect the results. 
Under Lipschitz and growth conditions on $b$ the moments and increment bounds
(\ref{momentbounds}, \ref{incrementbounds}) still hold true. In particular these 
imply that the solution $(t,x) \to X_t(x)$ is H\"{o}lder continuous. 
The extra term in the solution $ \int^t_0 \int G_{t-r}(x-u) b(X_r(u)) du dr$ is then 
continuously differentiable when $t>0, x \in \R$, and therefore will not affect the 
variations studied here.
This is also clear when the addition of
the drift induces an absolutely continuous change in the law of the solutions
(for example when $b$ is bounded and $\sigma(x) \geq \sigma_0 >0$).
Furthermore the results should go over for the equation (\ref{withdrift})
where $b(X_t(x))$ and $\sigma(X_t(x))$ are replaced by adapted fields
$b(t,x)$ and $\sigma(t,x)$ provided that they are known to be
H\"{o}lder continuous in $(t,x)$ and satisfy reasonable moment conditions. 

An example that does not fit into our Lipschitz assumptions
above is the equation (\ref{spde}) with $\sigma(x)=x^{\beta}$, for 
$\beta \in (0,1)$, and where one considers non-negative solutions. 
In particular the case $\beta=1/2$ occurs as the density of 
Dawson-Watanabe processes (see Konno and Shiga \cite{Konno+Shiga}). 
For these equations the moment
and increment bounds (\ref{momentbounds},\ref{incrementbounds}) still hold,
since the proofs rely only on the linear growth bounds on $\sigma$. 
The other place the Lipschitz assumption is needed is in Lemmas
\ref{keyspacelemma} and \ref{keytimelemma} where it is used
to get estimates on $E[|\sigma(X_t(x)) -\sigma(X_s(y))|]$.
For this we may use the inequality 
\[
|a^{\beta}-b^{\beta}|^2 \leq |a-b|^{2 \beta} \leq (1-\beta) \delta^{\gamma/2} 
+ \beta \delta^{-\gamma(1-\beta)/2\beta} |a-b|^2,
\]
valid for $a,b \geq 0$ and $\gamma,\delta>0$ (and derived from H\"{o}lder's inequality). 
Together with the increment bounds (\ref{incrementbounds}) one reaches, for 
$|t-s| \leq \delta^{\gamma}$ and $\lambda>0$,
\[ 
E \left[ |(X_t(x))^{\beta} - (X_s(y))^{\beta}|^2 \right]
\leq C(\beta, \lambda, T) \left( \delta^{\gamma \beta/2} + |x-y| \delta^{-\gamma(1-\beta)/2} \right)
e^{\lambda (|x|+|y|)}
\]
for all $x,y$ and $s,t \in [0,T]$. Using this in the proofs of Lemmas
\ref{keyspacelemma} and \ref{keytimelemma}, letting $t-t(\delta)=\delta^{\gamma}$ with 
$\gamma= 4/(2+\beta)$ and $\gamma=4/(4+\beta)$ respectively, 
the same variation results can be checked to hold. 

For the equation (\ref{spde}) over a finite interval
$x \in [a,b]$ with, say, periodic or Dirichlet or Neumann boundary conditions, 
the same variation results should hold on $(a,b) \times (0,\infty)$. 
This is intuitively clear since the small time
behaviour of the Green's function on the bounded interval remains the same as that 
of the heat kernel on the whole of $\R$. In the linear case, that is with constant 
coefficients, a 
simple way to deduce the results for bounded intervals from the whole space results 
is via a change of measure argument. In Mueller and Tribe \cite{Mueller+Tribe}
Corollary 4 there is an expression for the Radon Nikodym derivative of the law of $X_t(x)$
for a solution on $(a,b)$ with periodic boundary conditions with respect to a solution
on the whole space (with a periodically extended initial condition). Dropping to 
subsequence, if necessary, to obtain almost sure limits, the variation limits  
therefore still hold. (The proof that the expression is truly a change of measure is
incorrect in \cite{Mueller+Tribe}, but a simple corrected proof is 
given, as an errata, in the same journal). 


\begin{thebibliography}{99}
%
\bibitem{Barabasi+Stanley} A.~L.~Barabasi and H.~E.~Stanley, Fractal concepts in surface growth,
Cambridge University Press, Cambridge, 1995.
%
\bibitem{Giacomin} L.~Bertini and G.~Giambattista, Stochastic Burgers and KPZ equations 
from particle systems, {\em Comm. Math. Phys.\/} 183 (1997), no. 3, 571--607. 
%
\bibitem{Andref} R.~Carmona and S.~Molchanov, Parabolic Anderson problem and intermittency,
{\em Mem. Amer. Math. Soc.} 108 (1994), no. 518, viii+125 pp.
%
\bibitem{Konno+Shiga} N.~Konno and T.~Shiga, Stochastic partial differential
equations for some measure-valued diffusions, {\em Probab. Theory Related Fields\/} 79 (1988),
no. 2, 201--225.
%
\bibitem{Mueller+Tribe} C.~Mueller and R.~Tribe, Hitting properties of a 
random string, {\em Electron. J. Probab.\/} 7 (2002), no. 10, 29 pp. (electronic).
%
\bibitem{Pospisil} J.~Pospisil, Stochastic evolution equations driven by
fractional Brownian motion, PhD Thesis, University of West Bohemia, Plzen,
in preparation.
%
\bibitem{Shiga} T.~Shiga, Two contrasting properties of solutions for
 one-dimensional stochastic differential equations,
{\em Canad. J. Math.\/} 46 (1994), no. 2, 415--437.
%
\bibitem{Tribe} R.~Tribe, Large time behaviour of interface solutions to the heat equation with
Fisher--Wright white noise, {\em Probab. Theory Related Fields\/} 102 (1995), no. 3, 289--311.
%
\bibitem{Walsh} J.~B.~Walsh, An introduction to stochastic partial
differential equations, 
{\em \'Ecole d'\'et\'e de probabilit\'es de Saint-Flour, XIV---1984\/}, 265--439, 
Lecture notes in Math. 1180, Springer, Berlin, 1986.
%
\end{thebibliography}
\end{document}